\documentclass[ams,12pt]{article}
\usepackage{geometry}
\usepackage{amsmath,amssymb,amsfonts,amsthm}
\usepackage{txfonts}
\usepackage{todonotes}
\usepackage{tablists}
\usepackage{graphicx}
\usepackage{mathrsfs}
\usepackage{color}
\usepackage{float}
\usepackage{extarrows}
\usepackage{exscale}
\usepackage{titlesec}
\titleformat{\section}[block]{\large\center\sc}{\arabic{section}}{0.5em}{}[] 
\usepackage{relsize}
\usepackage{color}
\definecolor{teal}{RGB}{0,128,172}
\definecolor{pur}{RGB}{224,104,255}
\usepackage{hyperref}
\hypersetup{hypertex=true,
            colorlinks=true,
            linkcolor=teal,
            anchorcolor=blue,
            citecolor=pur}
\usepackage[titles]{tocloft} 
\usepackage{fancyhdr}
\theoremstyle{plain}
\newtheorem{theorem}{Theorem}[section]
\newtheorem{lemma}[theorem]{Lemma}

\newtheorem{remark}[theorem]{Remark}

\let\oldsection\section
\renewcommand\section{\setcounter{equation}{0}\oldsection}

\def\md{\mathrm{d}}
\def\be{\begin{equation}}
\def\ee{\end{equation}}
\def\bes{\begin{equation*}}
\def\ees{\end{equation*}}
\def\bs{\begin{split}}
\def\es{\end{split}}
\def\bali{\begin{aligned}}
\def\eali{\end{aligned}}
\newcommand{\pf}{\noindent {\bf Proof. \hspace{2mm}}}
\def\bR{{\mathbb R}}

\def\al{\alpha}
\def\e{\epsilon}

\def\la{\lambda}

\def\t{\tilde}
\def\th{\theta}

\def\Dl{\Delta}
\def\lt{\left}
\def\rt{\right}

\def\i{\infty}
\def\p{\partial}
\def\f{\frac}
\def\na{\nabla}

\def\O{\Omega}
\def\s{\sqrt}
\def\q{\quad}

\def\bl{\boldsymbol}

\def\mR{\mathbb{R}}

\def\fv{\mathfrak{v}}

\allowdisplaybreaks

\def\cd{\cdot}
\def\les{\lesssim}

\begin{document}

\title{\bf\Large  Existence of the planar stationary flow in the presence of interior sources and sinks in an exterior domain}

\author{\normalsize\sc Zijin Li and Xinghong Pan}

\date{}

\maketitle

\begin{abstract} In the paper, we consider the solvability of the two-dimensional Navier-Stokes equations in an exterior unit disk.  On the boundary of the disk, the tangential velocity is subject to the perturbation of a rotation, and the normal velocity is subject to the perturbation of an interior sources or sinks. At infinity, the flow stays at rest. We will construct a solution to such problem, whose principal part admits a critical decay $O(|x|^{-1})$. The result is related to an open problem raised by V. I. Yudovich in [{\it Eleven great problems of mathematical hydrodynamics}, Mosc. Math. J. 3 (2003), no. 2, 711--737], where Problem 2b states that: {\em Prove or disprove the global existence of stationary and periodic flows of a viscous incompressible fluid in the presence of interior sources and sinks.} Our result partially gives a positive answer to this open in the exterior disk for the case when the interior source or sink is a perturbation of the constant state.

\medskip



\end{abstract}

\tableofcontents

\section{Introduction}
We consider the 2D stationary Navier-Stokes flow in the planar exterior domain $\O:=\bR^2-B$, where $B=\{x\in\mathbb{R}^2\,:\,|x|<1\}$ is the unit disk.
\begin{equation}\label{NS}
\left\{\begin{array}{ll}
-\Delta \bl{u}+\bl{u} \cdot \nabla \bl{u}+\nabla p=\bl{f}, & x\in\O\,,\\[1mm]
\nabla \cdot \bl{u}=0, & x\in\O\,,\\[1mm]
\bl{u}=(\nu+g_r(\th))\bl{e}_r+(\mu+g_\th(\th))\bl{e}_\th, & x\in\p B\,,\\[1mm]
\bl{u}=0,&  |x|\rightarrow+\i.
\end{array}
\rt.
\end{equation}
 Here $\bl{u}$ is the unknown velocity of the fluid, while $p$ is the scalar pressure. $\bl{f}$ on the right hand is the external force. $\boldsymbol{n}$ and $\boldsymbol{\tau}$ are the unite outer normal vector and the tangential vector on the boundary $\p B$. $\bl{e}_r$ and $\bl{e}_\th$ are the orthogonal basis in the polar coordinates. $\mu$, $\nu\in\mathbb{R}$ are two constants, represents strengths of the rotation and flux. $\bl{g}:=g_r\bl{e}_r+g_\th\bl{e}_\th$ is the perturbation function defined on $\p B$ with suitable smoothness and small amplitude.  Without loss of generality, we assume that
 \be\label{zeromean}
\int^{2\pi}_0 g_r(\th)d\th=0\,,
 \ee
otherwise, we can rewrite
\[
\nu+g_r(\th)=\nu+\f{1}{2\pi}\int^{2\pi}_0 g_r(\th)d\th+g_r(\th)-\f{1}{2\pi}\int^{2\pi}_0 g_r(\th)d\th:=\t{\nu}+\t{g}_r(\th)\,,
\]
 where $\t{g}_r(\th):=g_r(\th)-\f{1}{2\pi}\int^{2\pi}_0 g_r(\th)d\th$ satisfies \eqref{zeromean}.

 The sign of $\nu$ determines that the flow has interior sources ($\nu>0$) or sink ($\nu<0$).

\begin{remark}\label{remflux}
Denote $B_r:=\{x\in\mathbb{R}^2\,:\,|x|<r\}$. The quantity $\Phi:=\int_{\p B_r} \bl{u}\cdot \boldsymbol{e}_r dS$ is defined as the flux flowing inside or outside, which is invariant with respect to $r$ by using the incompressible condition. Actually, by using the Gauss formula and the divergence-free condition, we have that for any $r\geq 1$,
\bes
\Phi:=\int_{\p B_r} \bl{u}\cdot \boldsymbol{e}_r dS= \int_{\p B} \bl{u}\cdot \boldsymbol{e}_r dS=2\pi\nu.
\ees
\end{remark}

\qed

Since our referenced domain is exterior to a disc, it is more convenient to reformulate system \eqref{NS} in polar coordinates. In polar coordinates $(r,\th)$, $x=r \cos \th,\ y=r \sin \th$, the polar orthogonal basis
\bes
\bl{e}_r=(\cos\th,\sin\th),\q \bl{e}_\th=(-\sin\th,\cos\th).
\ees
 In this polar coordinate system, we denote
\[
\bl{u}=u_r(r,\th)\bl{e_r}+u_\th(r,\th)\bl{e_\th},,
\]
and similarly as $\bl{f}$. In this way, we can rewrite \eqref{NS} as
\begin{equation}\label{PoNS}
\left\{\begin{split}
&-\left(\p_r^2+\f{1}{r}\p_r+\f{1}{r^2}\p_\th^2-\f{1}{r^2}\right)u_r+\f{2}{r^2}\p_\th u_\th+\p_rp+\left(u_r\p_r+\f{u_\th\p_\th}{r}\right)u_r-\f{u_\th^2}{r}=f_r\,,\\[1mm]
&-\left(\p_r^2+\f{1}{r}\p_r+\f{1}{r^2}\p_\th^2-\f{1}{r^2}\right)u_\th-\f{2}{r^2}\p_\th u_r+\f{1}{r}\p_\th p+\left(u_r\p_r+\f{u_\th\p_\th}{r}\right)u_\th+\f{u_ru_\th}{r}=f_\th\,,\\[1mm]
&\p_\th u_\th+\p_r(ru_r)=0,\\
&u_r(\th,1)=\nu +{g}_r(\th),\q u_\th(\th,1)=\mu+ g_\th(\th),\q (u_r,u_\th)\big|_{r\rightarrow+\i}=0.
\end{split}\right.
\end{equation}
where $g_r(\th)$ and ${g}_\th(\th)$ are two periodic function with respect to $\th$.

%

Before stating the main theorem in this paper, we need to define the functional spaces where we work. First for a periodic function $h(\th)\,:\,\mathbb{T}\to\mathbb{R}$, define its Fourier series by
$
h(\th)=\sum_{k\in\mathbb{Z}}h_k e^{ik\th},
$
where the sequence $\{h_k\}_{k\in\mathbb{Z}}$ is the set of its Fourier coefficients with $h_k=\f{1}{2\pi}\int^{2\pi}_0 h(\th) e^{-ik\th}d\th$. Also for a function $h(r):\, [1,+\i)\to\bR$, define its weighted $L^\i$ norm by
$\|h(r)\|_{L^\i_\zeta}:=\|r^\zeta h(r)\|_{L^\i([1,+\i))}$.

Now we define the first functional space which is for the solution of the reformualted system \eqref{PoNS}:
\be\label{BLa}
\mathcal{B}_\la =\left\{\bl{u}(\th,r):\,\mathrm{div}\,\bl{u}=0,\,\sum_{k\in\mathbb{Z},j\in\{r,\th\}}\left( (1+k^2)\|u_{j,k}(r)\|_{L^\i_{\la-2}}+(1+|k|)\|u'_{j,k}(r)\|_{L^\i_{\la-1}}+\|u''_{j,k}(r)\|_{L^\i_{\la}}\right)<\i\right\}\,,
\ee
where $\la>3$, and
\be\label{BLab}
\bar{\mathcal{B}}_\la=\left\{\bl{u}=\tilde{\bl{u}}+\f{\sigma}{r}\bl{e_\th}\,:\,\sigma\in\mR,\,\,\text{and}\,\,\tilde{\bl{u}}\in\mathcal{B}_\la\right\}\,.
\ee
Here the norm of ${\mathcal{B}}_\la$ and $\bar{\mathcal{B}}_\la$ are given by
\[
\begin{split}
&\|\bl{u}\|_{{\mathcal{B}}_\la}:=\sum_{k\in\mathbb{Z},j\in\{r,\th\}}\left( (1+k^2)\|u_{j,k}(r)\|_{L^\i_{\la-2}}+(1+|k|)\|u'_{j,k}(r)\|_{L^\i_{\la-1}}+\|u''_{j,k}(r)\|_{L^\i_{\la}}\right)\,,\\
&\|\bl{u}\|_{\bar{\mathcal{B}}_\la}:=|\sigma|+\|\tilde{\bl{u}}\|_{\mathcal{B}_\la}\,.
\end{split}
\]
For the convenience of the further expression, we denote
\[
\widetilde{\mathcal{B}}_\la=\left\{
\begin{split}
\mathcal{B}_\la,\q\text{for }\q\nu<-2\,;\\
\bar{\mathcal{B}}_\la,\q\text{for }\q\nu\geq-2\,.\\
\end{split}
\right.
\]
Meanwhile, we define the following two spaces, which are for the external force and the boundary value, respectively:
\[
\begin{split}
\mathcal{E}_\la &=\left\{\bl{f}(\th,r)\,:\,\mathbb{T}\times[1,\i)\to\mathbb{R}^2\,\big|\,\sum_{k\in\mathbb{Z},j\in\{r,\th\}}\|f_{j,k}(r)\|_{L^\i_\la}<\i\right\}\,;\\
\mathcal{V}&=\left\{\bl{g}(\th)\,:\,\mathbb{T}\to\mathbb{R}^2\,\big|\,\sum_{k\in\mathbb{Z},j\in\{r,\th\}}(1+k^2)|g_{j,k}|<\i\right\}\,.
\end{split}
\]
And their norms are defined by
\[
\begin{split}
\|\bl{f}\|_{\mathcal{E}_\la} &:=\sum_{k\in\mathbb{Z},j\in\{r,\th\}}\|f_{j,k}(r)\|_{L^\i_\la}\,;\\
\|\bl{g}\|_{\mathcal{V}}&:=\sum_{k\in\mathbb{Z},j\in\{r,\th\}}(1+k^2)|g_{j,k}|\,.
\end{split}
\]
Below is our main theorem:
\begin{theorem}\label{Main}
 Assume that
 \[
 \nu\leq-\f{3}{2}\,;
 \]
 or alternatively
 \[
 \nu>-\f{3}{2}\q\&\q|\mu|>\s{2\nu^3+19\nu^2+56\nu+48}\,.
 \]
There exists $\e>0$, being sufficiently small and depending on $\mu$ and $\nu$, such that if
\[
\|\bl{f}\|_{\mathcal{E}_\la}+\|\bl{g}\|_{\mathcal{V}}<\e\,,
\]
then problem \eqref{PoNS} has a unique solution $\bl{u}$ such that
\[
\bl{u}=\f{\nu}{ r}\bl{e}_r+\f{{\mu}}{r}\bl{e}_\th+\bl{v}\,,
\]
where $\bl{v}\in\widetilde{\mathcal{B}}_\la\,$ that satisfies
\[
\|\bl{v}\|_{\widetilde{\mathcal{B}}_\la }\leq C\left(\|\bl{f}\|_{\mathcal{E}_\la }+\|\bl{g}\|_{\mathcal{V}}\right)\leq C\e\,
\]
for some constants $\la=3^+$ and $C>0$, depending on $\mu$ and $\nu$. Here $3^+$ denotes a constant which is larger but close to $3$.

\end{theorem}

\qed

\begin{remark}
  Yudovich in Problem 2b of \cite{Yudovich:2003MMJ} states that ``Prove or disprove the global existence of stationary and periodic flows of a viscous incompressible fluid in the presence of interior sources and sinks."  Our result partially gives a positive answer to this problem in the exterior disc for the case when the interior source (Corresponding to $\Phi>0$) or sink (Corresponding to $\Phi<0$) is a perturbation of the constant state.
\end{remark}

\begin{remark}
Here we mention two papers. This first one is Higaki \cite{Higaki2023}, where a similar problem as \eqref{NS} is considered, but there the author only considered the case $g_\tau=g_n\equiv 0$ and the flow has an interior sink ($\nu<-2$). The second one is Hillairet and Wittwer \cite{Hillairet2013}, where the case $\nu=0$, $\int^{2\pi}_0 g_r(\th)d\th=0$ and $\bl{f}\equiv0$ is considered, which means the flow have neither interior sources or sinks, nor external forces.
\end{remark}

\qed

\begin{remark}
When the current paper has been finished, we found a similar result has been obtained in \cite{GH2024} independently. Although there are some overlaps between the results there and in this paper, the framework of analysis is different in many aspects.
\end{remark}

\qed

\begin{remark}
Actually, from the proof of Theorem \ref{Main} in the case of $\nu<-2$, we can choose $\t{\mu}$, which is different but sufficiently close to $\mu$, such that
  \[
\t{\bl{u}}=\f{\nu}{ r}\bl{e}_r+\f{\t{\mu}}{r}\bl{e}_\th+\bl{v},\q \bl{v}\in \mathcal{B}_\la\,,
\]
 is still a solution of system \eqref{NS}. Since $\bl{v}\in\mathcal{B}_\la$ is a subcritically decayed space, then we obtain the non-uniqueness of solutions to system \eqref{NS} in the case $\nu<-2$. Here we thank Prof. Guo and Prof. Hillairet, who informed us a similar result in their earlier paper \cite{GH2024}.  However, whether there is non-uniqueness for the problem \eqref{NS} with $\nu\geq-2$ is still unknown to the authors.
\end{remark}

\qed

The existence problem for system \eqref{NS} is closely connected to the existence problem of the 2D exterior-domain problem, which states that to fine a solution to the following problem
\begin{equation}\label{ns2}
\left \{
\begin {array}{ll}
\bl{u}\cdot\na \bl{u}+\na p-\Dl \bl{u}=0, & \text{in }\bR^2-D\,,\\ [5pt]
\na\cdot\bl{u}=0, & \text{in }\bR^2-D\,,\\
\bl{u}\big|_{\p D}=\bl{a}^\ast,\\
\bl{u}\big|_{|x|\rightarrow+\i}=\eta\bl{e}_1,\\
\end{array}
\right.
\end{equation}
where $D$ is a smooth bounded domain and $\bl{a}^\ast$ is a smooth function defined on $\p D$. The constant is to distinguish the case of a flow around $D$ ($\beta=0$) and a flow past $D$ ($\eta\neq 0$).

 The existence of the 2D exterior domain problem was paid attention to since the Stokes paradox, which states that when considering the linear Stokes equation of \eqref{ns2}, this is no solution. For the Navier-Stokes system, with general $\bl{a}^\ast$ and $\eta$, such an existence problem of \eqref{ns2} was listed by Yudovich in \cite{Yudovich:2003MMJ} as one of the ``Eleven Great Problems in Mathematical Hydrodynamics" (Problem 2), which was initially studied by Leray in \cite{Leray:1933JMPA} by using the {\em invading domains method}. By using Leray's method, a $D$-solution (the solution have finite Dirichlet integration) satisfying \eqref{ns2}$_{1,2,3}$ and no flux condition $\int_{\p D} \bl{a}^\ast \cdot \bl{n} dS=0$ can be obtained in \cite{KorobkovPR:2020JDE}. However, whether this $D$-solution satisfies \eqref{ns2}$_4$ is unknown. Also, if the flux of the flow is non-zero, whether a $D$-solution is existed to satisfy \eqref{ns2}$_{1,2,3}$ is still not clear. The main difficulties of the existence of the 2D exterior domain lie in the following two factors: The lack of Sobolev embedding in two dimensions and the logarithmic growth of the Green tensor for the 2D Stokes system. Although the above difficulties, Finn and Smith in \cite{FinnS:1967ARMA} gave an existence result for system \eqref{ns2} in  the case that $\eta$ and $|\bl{a}^\ast-\eta\bl{e}_1|$ is small with the help of iteration techniques. Whether the Finn-Smith solution is a $D$-solution stays unknown. Recently Korobkov-Ren \cite{KorobkovR:2021ARMA, KorobkovR:2022JMPA} shows existence and uniqueness of $D$-solutions in the case that $a^\ast=0$ and $\eta$ is small. See recent advances on this topic in \cite{KorobkovPR:2019ARMA,KorobkovPR:2021ADV} and references therein. Our main result in Theorem \ref{Main} gives a solution to the 2D exterior domain problem \eqref{NS}, which can have arbitrary flux. Also the constructed solution is a $D$-solution.

When the external force $\bl{f}$ and the perturbation $(g_\tau,g_n)$ are trivial, there is an explicit solution $\f{\nu}{r} \bl{e}_r+\f{\mu}{r}\bl{e}_\th$ to system \eqref{NS},
which is invariant under the natural scaling of the Navier-Stokes equations: $\bl{u}^{\al}(x):=\al\bl{u}(\al x)$. A scaling-invariant solution is called scale-critical and represents the balance between the nonlinear and linear parts of the equations. Given this nature, perturbation around a scaling-invariant solution maybe complicated based on the scale of the perturbation. It is expected that the problem is well-posed if the perturbation is subcritical. Hillairet and Wittwer \cite{Hillairet2013} consider the perturbation of system \eqref{NS} around $\mu\bl{e}_\th$ in an exterior disk. Also the flow is assumed to be  zero flux and zero external force.  They show that when $|\mu|>\s{48}$, the linearized equations for the vorticity fall into the subcritical category. Then, iteration to the nonlinear problem with subcritical nonlinearity can be closed to produce a subcritical vorticity. they show the existence of solutions in the form of $\bl{u}=\mu\bl{e}_\th+o(|x|^{-1})$ when $|x|\rightarrow+\i$. See some related results in Higaki \cite{HigakiMN:2018ARMA,GallagherHM:2019MATHNA} for the flow with zero flux, non-perturbed boundary condition and non-zero external force. Recently Higaki \cite{Higaki2023} consider the external force perturbation effect of system \eqref{NS} with $\nu<-2$, in which the boundary condition is not perturbed. (corresponding to $g_\tau=g_n=0$ in \eqref{NS}$_3$). A similar result as Theorem \ref{Main} was obtained. Our result can be viewed as an improvement to the above-mentioned results. The main result in  Theorem \ref{Main} can be explained as the scale-critical flow  $\f{\nu}{r}\bl{e}_r+\f{\mu}{r}\bl{e}_\th$ can produce a stabilizing effect to the spatial decay when $\mu,\nu$ satisfy some suitable constraints.
\\[2mm]

{\noindent\bf Strategy of Proof to the main result}

The strategy of proving Theorem \ref{Main} are the following. First we construct the solution $(\bl{u},p)$ of system \ref{NS} in the form of
\bes
\bl{u}=\f{\nu}{r}\bl{e}_r+\f{{\mu}}{r}\bl{e}_\th+\bl{v}.
\ees
The error $\bl{v}$ is understood as the perturbation from of $\f{\nu}{r}\bl{e}_r+\f{{\mu}}{r}\bl{e}_\th$  in response to the external force $\bl{f}$ and the boundary condition. Then using the relation
\bes
\bl{u}\cdot\na \bl{u}=\bl{u}^{\bot} \na\times\bl{u}+\na\lt(\f{|\bl{u}|^2}{2}\rt),\q \na\times\bl{u}:=\p_1u_2-\p_2 u_1,\q \bl{u}^{\bot}:=(-u_2,u_1),
\ees
we see that $\bl{v}$ satisfies the following error system
\begin{align}\label{nserror}
\lt\{
\begin{array}{ll}
-\Dl \bl{v}+\lt(\f{\nu}{r}\bl{e}_r+\f{{\mu}}{r}\bl{e}_\th \rt)^{\bot}\na\times\bl{v}+\na \pi=-\bl{v}^{\bot}\na\times\bl{v}+\bl{f},& \text{in } \O,\\[1mm]
\na\cdot\bl{v}=0,& \text{in } \O,\\[1mm]
v_\th= g_\th(\th),\q v_r={g}_r(\th)\,,& \text{on } \p B, \\[1mm]
\bl{v}=0, & \text{as } r\rightarrow+\i,
\end{array}
\rt.
\end{align}
where $\na \pi=\na\lt(p+\f{|\bl{u}|^2}{2}\rt)$.

Next, we show the existence of the solution to the problem \eqref{nserror}, which depends on the linear structure as follows.
\begin{align}\label{nserrorlinear}
\lt\{
\begin{array}{ll}
-\Dl \bl{\fv}+\lt(\f{\nu}{r}\bl{e}_r+\f{{\mu}}{r}\bl{e}_\th \rt)^{\bot}\na\times\bl{\fv}+\na q=\bl{f},& \text{in } \O\,,\\[1mm]
\na\cdot\bl{\fv}=0,& \text{in } \O\,,\\[1mm]
\fv_\th= g_\th(\th),\q \fv_r={g}_r(\th),& \text{on } \p B\,, \\[1mm]
\bl{\fv}=0, & \text{as } r\rightarrow+\i\,.
\end{array}
\rt.
\end{align}
One can study the linearized system \eqref{nserrorlinear} in each Fourier mode. System \eqref{nserrorlinear} will be solved for zero mode and non-zero mode separately.

For the zero mode, $\bl{v}_0=v_{r,0}\bl{e}_r+v_{\th,0}\bl{e}_r$, it is  not hard to deduce that $v_{r,0}=0$ by using the boundary condition \eqref{nserrorlinear}$_3$ and the incompressible condition. While $v_{\th,0}$ satisfies the following ODE
\be\label{homzero}
-\left(\f{\md^2}{\md r^2}+\f{1-\nu}{r}\f{\md}{\md r}-\f{1+\nu}{r^2}\right)\fv_{\th,0}=f_{\th,0}\,,
\ee
the two linearly independently fundamental solutions for the homogeneous equation of \eqref{homzero} are: $r^{\nu+1}$ and $r^{-1}$ (or $r^{-1}$ and $r\log r$ for $\nu=-2$). Thus for $\nu<-2$, if $f_{\th,0}$ decays subcritically, we can obtain a solution of \eqref{homzero} which decay subcritically and satisfying the prescribed boundary condition $v_{\th,0}=g_{\th,0}$. However, when $\nu\geq2$, the fundamental solutions $r^{\nu+1}$ or $r\log r$ is supercritical. Although we can construct a subcritically decayed solution $\t{v}_{\theta,0}$ of \eqref{homzero} when $f_{\th,0}$ decays fast enough at spatial infinity, it may not satisfy the boundary condition. In this situation we need to correct this solution with a critical decay term $\f{\sigma}{r}$ such that it can fulfill the boundary condition. That is why there is an extra $\sigma$ in the definition of functional space $\bar{\mathcal{B}}_\la$ when considering the case of $\nu\geq-2$.

For the non-zero mode of \eqref{nserrorlinear}, we derive the representation formula for the $k-$mode ($k\neq0$) of the velocity by using the stream and vorticity functions. In this  way we overcome the difficulties caused by the pressure. By solving the vorticity equation, we can recover the velocity by the Biot-Savart law. In order to using contract mapping to obtain a solution of the nonlinear problem, the stream and vorticity functions need to decay subcritically, which require that $\mu$ and $\nu$ satisfy the constraints in Theorem \ref{Main}.
%

%

Our paper is organized as follows. In section \ref{sec2}, we formulate the linearized system \eqref{nserrorlinear} in polar coordinates and deduce its each Fourier mode. The zero mode is direct from the equation, while the non-zero mode is recovered from the vorticity equation. In section \ref{sec3}, we will solve each Fourier mode in subcritically decayed function space. Then by using contract mapping, we show the existence of solutions to the nonlinear system in section \ref{sec4}.

Throughout the paper, $C_{a,b,...}$ denotes a positive constant depending on $a,\,b,\,...$, which may be different from line to line. $A\lesssim B$ means $A\leq CB$.

\section{Formulation of the linear system in Polar coordinates and Fourier mode}\label{sec2}

First, we reformulate the nonlinear system \eqref{nserror} in polar coordinates. Denoting $\bl{v}:=v_r\bl{e}_r+v_\th\bl{e}_\th$, we can deduce
\begin{small}
\[
\left\{\begin{split}
&-\left(\p_r^2+\f{1-\nu}{r}\p_r+\f{1}{r^2}\p_\th^2-\f{1}{r^2}\right)v_r+\f{\mu}{r^2}\p_\th v_r+\f{2}{r^2}\p_\th v_\th-\f{2\mu}{r^2}v_\th+\p_r\pi=-\left({v}_r\p_r+\f{{v}_\th}{r}\p_\th\right){v}_r+\f{{v}_\th^2}{r}+f_r\,,\\[1mm]
&-\left(\p_r^2+\f{1-\nu}{r}\p_r+\f{1}{r^2}\p_\th^2-\f{1+\nu}{r^2}\right)v_\th+\f{\mu}{r^2}\p_\th v_\th-\f{2}{r^2}\p_\th v_r+\f{1}{r}\p_\th \pi=-\left({v}_r\p_r+\f{{v}_\th}{r}\p_\th\right){v}_\th-\f{{v}_r{v}_\th}{r}+f_\th\,,\\[1mm]
&\p_\th v_\th+\p_r(rv_r)=0,\\
&v_\th\big|_{r=1}= g_\th(\th),\q v_r\big|_{r=1}={g}_r(\th), \q \bl{v}\big|_{r\rightarrow+\i}=0,
\end{split}\right.
\]
\end{small}
Then the related  linearized system \eqref{nserrorlinear} in polar coordinates, after denoting $\bl{\fv}:=\fv_r\bl{e}_r+\fv_\th\bl{e}_\th$, reads
\be\label{PoNSvL}
\left\{\begin{split}
&-\left(\p_r^2+\f{1-{\nu}}{r}\p_r+\f{1}{r^2}\p_\th^2-\f{1}{r^2}\right)\fv_r+\f{\mu}{r^2}\p_\th \fv_r+\f{2}{r^2}\p_\th \fv_\th-\f{2\mu}{r^2}\fv_\th+\p_r\pi=f_r\,,\\[1mm]
&-\left(\p_r^2+\f{1-{\nu}}{r}\p_r+\f{1}{r^2}\p_\th^2-\f{1+{\nu}}{r^2}\right)\fv_\th+\f{\mu}{r^2}\p_\th \fv_\th-\f{2}{r^2}\p_\th \fv_r+\f{\p_\th}{r}\pi=f_\th\,,\\[1mm]
&\p_\th \fv_\th+\p_r(r\fv_r)=0\,,\\
&\fv_r\big|_{r=1}={g}_r(\th),\q\fv_\th\big|_{r=1}= g_\th(\th),\q \bl{\fv}\big|_{r\rightarrow+\i}=0\,.
\end{split}\right.
\ee
Applying the Fourier series technique, i.e.
\[
\fv_r(r,\th):=\sum_{k\in\mathbb{Z}}\fv_{r,k}(r)e^{ik\th}\,,\q\text{where}\q\fv_{r,k}(r)=\f{1}{2\pi}\int_0^{2\pi}\fv_r(r,\th)e^{-ik\th}\md\th\,,
\]
and similarly for $\fv_\th$ and the other functions, we rewrite \eqref{PoNSvL} in each $k$-mode for any $k\in\mathbb{Z}$ :
\be\label{PoNSvLF}
\left\{\begin{split}
&-\left(\f{\md^2}{\md r^2}+\f{1-{\nu}}{r}\f{\md}{\md r}-\f{1+k^2}{r^2}\right)\fv_{r,k}+\f{ik\mu}{r^2}\fv_{r,k}+\f{2ik}{r^2} \fv_{\th,k}-\f{2\mu}{r^2}\fv_{\th,k}+\p_rq_k=f_{r,k}\,,\\[1mm]
&-\left(\f{\md^2}{\md r^2}+\f{1-{\nu}}{r}\f{\md}{\md r}-\f{1+{\nu}+k^2}{r^2}\right)\fv_{\th,k}+\f{ik\mu}{r^2}\fv_{\th,k}-\f{2ik}{r^2}\fv_{r,k}+\f{ik}{r}q_k=f_{\th,k}\,,\\[1mm]
&ik\fv_{\th,k}+\p_r(r\fv_{r,k})=0\,,\\[1mm]
&\fv_{r,k}(1)=g_{r,k}\,,\q \fv_{\th,k}(1)=g_{\th,k}\,\q \fv_{r,k}(+\i)=\fv_{\th,k}(+\i)=0.
\end{split}\right.
\ee

\subsection*{The Fourier zero mode}
From \eqref{PoNSvLF}$_{3,4}$, we see that  $\p_r(r\fv_{r,0})=0$ and $\fv_{r,0}(1)=\fv_{r,0}(+\i)=0$, which indicates that $\fv_{r,0}\equiv 0$. From \eqref{PoNSvLF}$_{2,4}$, we see that $v_{\th,0}$ satisfies

\be\label{PoNSvLF0}
\left\{\begin{split}
&-\left(\f{\md^2}{\md r^2}+\f{1-\nu}{r}\f{\md}{\md r}-\f{1+\nu}{r^2}\right)\fv_{\th,0}=f_{\th,0}\,,\q\q\q\q\q\text{for}\q r\in(1,\i)\,,\\[1mm]
&\fv_{\th,0}(1)=g_{\th,0}\,\q \fv_{\th,0}(+\i)=0.
\end{split}\right.
\ee
\subsection*{The case of $\bl{k\in\mathbb{Z}-\{0\}}$}
To overcome the difficulty caused by the pressure term in the case $k\neq0$, we introduce the stream function and vorticity to the linearized system \eqref{nserrorlinear}. Since the related linear velocity field $\bl{\mathfrak{v}}:=\fv_r\bl{e_r}+\fv_\th\bl{e_\th}$ is divergence-free and the zero mode of $\fv_r$ is zero, there exists a periodic stream function $\phi$ in $\th$ variable such that
\be\label{BS}
\fv_r=\f{\p_\th\phi}{r},\q\text{and}\q \fv_\th=-\p_r\phi\,,
\ee
and define the vorticity $w$ by
\be\label{yphi}
w:=\f{1}{r}\p_r(r\fv_\th)-\f{1}{r}\p_\th \fv_r\,,
\ee
which satisfies
\be\label{Ephi}
-\Dl \phi=w.
\ee
Using \eqref{PoNSvL}, we derive $w$ satisfies
\be\label{xephi}
-\left(\p_r^2+\f{1-\nu}{r}\p_r+\f{1}{r^2}\p_\th^2-\f{\mu}{r^2}\p_\th\right)w=F:=\f{1}{r}\p_r(rf_\th)-\f{1}{r}\p_\th f_r\,.
\ee

Splitting \eqref{Ephi} and \eqref{xephi} in each independent Fourier mode, we arrive at
\begin{align}\label{streamvor}
\lt\{
\bali
&-\left(\f{\md^2}{\md r^2}+\f{1}{r}\f{\md}{\md r}-\f{k^2}{r^2}\right)\phi_k=w_k\,,\\
&-\left(\f{\md^2}{\md r^2}+\f{1-\nu}{r}\f{\md}{\md r}-\f{i\mu k+k^2}{r^2}\right)w_k=F_k:=\f{1}{r}\f{d}{dr}(rf_{\th,k})-\f{ik}{r}f_{r,k}\,,\\
&\phi_k(1)=\bar{\phi}_k,\q w_k(1)=\bar{w}_k\,,\\
& \phi_k(+\i)=w_k(+\i)=0\,,
\eali
\rt.\q\q\q \text{for}\q k\in\mathbb{Z}-\{0\}\,.
\end{align}
The boundary conditions for $\phi_{k}$ and $w_k$ on $r=1$ are subtle, which will be specialized later. The vanishing boundary condition at infinity for the stream function is to ensure that the solution decays subcritically. From \eqref{yphi}, applying the Fourier series technique, we find in each mode
\[
w_k(r):=\f{1}{r}\f{\md}{\md r}(r\fv_{\th,k}(r))-\f{ik}{r}\fv_{r,k}(r)\,.
\]
So \eqref{streamvor}$_2$  can also be obtained from \eqref{PoNSvLF} by cancelling $q_k$.  After solving \eqref{streamvor}, then the non-zero mode of $v_{r,k}$ and $v_{\th,k}$ is recovered from \eqref{BS} by
\be\label{stream1}
\fv_{r,k}=\f{ik\phi_k}{r},\q \fv_{\th,k}=-\phi'_k\,.
\ee

\section{Solvability of the Fourier mode in subcritically decayed spaces}\label{sec3}

In order to close the nonlinear decay estimates, the solutions we constructed need to decay subcritically, namely, the decay order of the stream function $\phi$, $\bl{\fv}$ and vorticity $w$ with respect to $r$ at spacial infinity are $-0^{+}$, $-1^{+}$ and $-2^{+}$. In our following analysis, we assume that the force $\bl{f}=f_r\bl{e}_r+f_\th\bl{e}_\th$ stays in $\mathcal{E}_\la$ and the boundary condition $\bl{g}=g_r\bl{e}_r+g_\th\bl{e}_\th$ stays in $\mathcal{V}$ for some $\la=3^+$, which will be determined later on.

\subsection*{Solvability for the zero mode}

\begin{lemma}\label{lemzeromode}
The linearized Fourier zero mode equation \eqref{PoNSvLF0} have a solution $\fv_{\th,0}$ satisfying
\[
\fv_{\th,0}=\t{\fv}_{\th,0}+\f{\sigma}{r}\chi_{\nu\geq-2},
\]
where $\t{\fv}_{\th,0}$ is a subcritically decayed solution of \eqref{PoNSvLF0}$_1$ and $\chi_{\nu\geq-2}$ is the characteristic function on $\{\nu\geq -2\}$. Also there exists a $\la_0>3$ such that for any $\la\in(3,\la_0).$
\[
\|\t{\fv}_{\th,0}\|_{L^\i_{\la-2}}+\|\t{\fv}'_{\th,0}\|_{L^\i_{\la-1}}+\|\t{\fv}''_{\th,0}\|_{L^\i_{\la}}+|\sigma|\leq C_{\nu,\la}\left(\|f_{\th,0}\|_{L^\i_{\la}}+|g_{\th,0}|\right).
\]
\end{lemma}
\pf Direct calculation shows the homogeneous equation of \eqref{PoNSvLF0}$_2$
\[
-\frac{\mathrm{d}^2 {\fv}_{\theta, 0}}{\mathrm{~d} r^2}-\frac{1-\nu}{r} \frac{\mathrm{~d} {\fv}_{\theta, 0}}{\mathrm{~d} r}+\frac{1+\nu}{r^2}{\fv}_{\theta, 0}=0\,
\]
has two linearly independent solutions: $r^{\nu+1}$ and $r^{-1}$ (or $r^{-1}$ and $r\log r$ for $\nu=-2$).

{\bf\noindent Case 1, $\bl{\nu<-2}$}

Direct solve the ODE \eqref{PoNSvLF0}$_2$ with boundary condition $\fv_{\th,0}(1)=g_{\th,0}$, which has the following solution decaying subcritically at spacial innifity
\be\label{v0<}
\begin{aligned}
\fv_{\theta, 0}(r)=&-\frac{1}{\nu+2}\left\{r^{\nu+1} \int_1^r s^{-\nu} f_{\theta, 0}(s) \mathrm{d} s+r^{-1} \int_r^{\infty} s^2 f_{\theta, 0}(s) \mathrm{d} s\right\}\\
                 &+\left(g_{\th,0}+\frac{1}{\nu+2}\int_1^{\infty} s^2 f_{\theta, 0}(s) \mathrm{d} s\right)r^{\nu+1},
\end{aligned}
\ee
which behaves as $O(r^{-1^{+}})$ as $r\to\i$. Since
$
|{f}_{\th,0}(r)|\leq \|f_{\th,0}\|_{L^\i_{\la}}r^{-\la}
$
for $3<\la<1-\nu:=\la_0$, by \eqref{v0<}, one derives
\be\label{Evth0}
\begin{split}
|r^{\la-2}\fv_{\th,0}(r)|&\leq-\f{\|f_{\th,0}\|_{L^\i_{\la}}}{\nu+2}\left(r^{\la+\nu-1}\int_1^rs^{-\nu-\la}\mathrm{d} s+r^{\la-3} \int_r^{\infty} s^{2-\la}\mathrm{d} s+r^{\la+\nu-1}\int_1^\i s^{2-\la}\mathrm{d}s\right)+|g_{\th,0}|r^{\la+\nu-1}\\[2mm]
&=-\f{\|f_{\th,0}\|_{L^\i_{\la}}}{\nu+2}\left(\f{1}{-\nu-\la+1}(1-r^{\nu+\la-1})+\f{1}{\la-3}(1+r^{\nu+2})\right)+|g_{\th,0}|r^{\nu+\la-1}\\
&\leq -\f{1}{\nu+2}\left(\f{1}{1-\nu-\la}+\f{2}{\la-3}\right)\|f_{\th,0}\|_{L^\i_{\la}}+|g_{\th,0}|\,.
\end{split}
\ee
Meanwhile, direct calculation shows
\[
\begin{split}
\f{\md}{\md r}\fv_{\theta, 0}(r)=&-\f{1}{\nu+2}\left\{(\nu+1)r^\nu\int_1^r s^{-\nu}f_{\th,0}(s)\md s-r^{-2}\int_r^\i s^2f_{\th,0}(s)\md s\right\}\\
&+(\nu+1)\left(g_{\th,0}+\frac{1}{\nu+2}\int_1^{\infty} s^2 f_{\theta, 0}(s) \mathrm{d} s\right)r^{\nu}\,.
\end{split}
\]
And similarly as one derives \eqref{Evth0} and using \eqref{PoNSvLF0}, we have
\be\label{vthe0first}
|r^{\la-1}\fv'_{\th,0}(r)|+|r^{\la}{\fv}''_{\th,0}(r)|\leq C_{\nu,\la}\left(\|f_{\th,0}\|_{L^\i_{\la}}+|g_{\th,0}|\right)\,.
\ee
Combining \eqref{Evth0} and \eqref{vthe0first}, we see that there exists a constant $C_{\la,\nu}$ such that
\[
\|{\fv}''_{\th,0}\|_{L^\i_{\la}}+\|\fv'_{\th,0}\|_{L^\i_{\la-1}}+\|\fv_{\th,0}\|_{L^\i_{\la-2}}\leq C_{\nu,\la}\left(\|f_{\th,0}\|_{L^\i_{\la}}+|g_{\th,0}|\right).
\]

{\bf\noindent Case 2, $\bl{\nu\geq-2}$}

However, if $\nu\geq-2$, the function $r^{\nu+1}$ (or $\log r$ when $\nu=-2$) decays slower than the prescribed request, one cannot solve the ODE of $\fv_{\th,0}$ to get a solution decays faster than $r^{-1}$ for any given $g_{\th,0}$. Instead, one has the following exact solution of \eqref{PoNSvLF0}$_2$
\be\label{v0>}
\fv_{\th,0}(r)=-\f{1}{r}\int_r^\i s^{\nu+1}\int_s^\i t^{-\nu}f_{\th,0}(t)\mathrm{d}t\mathrm{d}s+\f{\sigma}{r}:=\t{\fv}_{\th,0}(r)+\f{\sigma}{r},
\ee
where $\sigma$ satisfies
\[
\sigma-\int_1^\i s^{\nu+1}\int_s^\i t^{-\nu}f_{\th,0}(t)\mathrm{d}t\mathrm{d}s=g_{\th,0}\,.
\]
Owing to $f_{\th,0}\in L^\i_{\la}$ (Set $\la\in (3,\la_0)$, where $\la_0=3.01$.), it is clear that
\[
|\t{\fv}_{\th,0}(r)|\les r^{-\la+2}\|f_{\th,0}\|_{L^\i_\la}\,,
\]
which indicates $\t{\fv}_{\th,0}(r)$ decays faster than $r^{-1}$ as $r\to\i$. And it satisfies
\be\label{vth0sec}
\begin{split}
|r^{\la-2}\t{\fv}_{\th,0}(r)|&\leq r^{\la-3}\|f_{\th,0}\|_{L^\i_{\la}}\int_r^\i s^{\nu+1}\int_s^\i t^{-\nu-\la}\mathrm{d}t\mathrm{d}s=\f{1}{(\nu+\la-1)(\la-3)}\|f_{\th,0}\|_{L^\i_{\la}}\,.\\
\end{split}
\ee
Meanwhile, since
\[
\f{\md}{\md r}\t{\fv}_{\th,0}(r)=\f{1}{r^2}\int_r^\i s^{\nu+1}\int_s^\i t^{-\nu}f_{\th,0}(t)\mathrm{d}t\mathrm{d}s+r^\nu\int_r^\i t^{-\nu}f_{\th,0}(t)\mathrm{d}t\,,
\]
it is clear that
\be\label{vth0thi}
|r^{\la-1}\t{\fv}'_{\th,0}(r)|\leq C_{\nu,\la}\|f_{\th,0}\|_{L^\i_{\la}}\,.
\ee
Combining \eqref{v0>}, \eqref{vth0sec} and \eqref{vth0thi}, we can obtain that
\bes
\fv_{\th,0}(r)=\t{\fv}_{\th,0}(r)+\f{\sigma}{r},
\ees
where
\bes
\|\t{\fv}''_{\th,0}\|_{L^\i_{\la}}+\|\t{\fv}'_{\th,0}\|_{L^\i_{\la-1}}+\|\t{\fv}_{\th,0}\|_{L^\i_{\la-2}}+|\sigma|\leq C_{\nu,\la}\left(\|f_{\th,0}\|_{L^\i_{\la}}+|g_{\th,0}|\right).
\ees

\qed



\subsection*{Solvability for the non-zero mode}

\begin{lemma}\label{lemnonzeromode}
Under the assumption for $\nu$ and $\mu$ in Theorem \ref{Main}, there exists a $\la_1>3$ such that for $\la\in(3,\la_1)$. The linearized Fourier non-zero mode equation \eqref{streamvor} with \eqref{stream1} have a solution $\bl{\fv}_{k}$ satisfying
\bes
k^2\|\bl{\fv}_{k}\|_{L^\i_{\la-2}}+|k|\cd\|\bl{\fv}'_{k}\|_{L^\i_{\la-1}}+\|\bl{\fv}''_{k}\|_{L^\i_{\la}}\leq C_{\mu,\nu}k^2|\bl{g}_k|+ C_{\mu,\nu}\|\bl{f}_k\|_{L_\lambda^{\infty}}\,,\q\text{for all}\q k\in\mathbb{Z}-\{0\}\,.
\ees
\end{lemma}

\qed

\pf For the second one of \eqref{streamvor}, it is an Eulerian ODE whose homogeneous equation has two linearly independent solutions
\[
w_\pm(r)=r^{\xi_k^{\pm}},\q  \text{for } \xi_k^{\pm}=\f{\nu\pm\sqrt{\nu^2+4(k^2+i\mu k)}}{2}.
\]
 Direct calculation shows
\begin{align}
&\mathrm{Re\,}\xi_k^\pm=\f{\nu}{2}\pm\f{1}{2\sqrt{2}}\left[\left((\nu^2+4k^2)^2+(4\mu k)^2\right)^{1/2}+\nu^2+4k^2\right]^{1/2}, \label{xifirst}\\
&{\small |\xi_k^{\pm}|=\f{1}{2}\left(\nu^2\pm\sqrt{2}\nu\left[\left((\nu^2+4k^2)^2+(4\mu k)^2\right)^{1/2}+\nu^2+4k^2\right]^{1/2}+\left((\nu^2+4k^2)^2+(4\mu k)^2\right)^{1/2}\right)^{1/2}.}\label{xisecond}
\end{align}

Recalling that we intend to solve a subcritically decayed solution such that insist $w_k(r)=o(r^{-2})$ as $r\to\i$. Clearly $\mathrm{Re\,}\xi_k^-$ is increasing with $|k|$, so we only require ${\mathrm{Re}\,\xi_1^{-}<-2}$. In this situation, by direct calculate \eqref{xifirst}, we can obtain that
\be
\nu<-\f{3}{2},
\ee
or
\[
\nu\geq -\f{3}{2},\q\text{and}\q  |\mu|>\s{2\nu^3+19\nu^2+56\nu+48}\,,
\]
which is the requirement in Theorem \ref{Main}. In this case, the subcritically decayed solution of \eqref{streamvor}$_{2,3,4}$ is given by
\begin{align}
w_k(r)=&{\bar{w}_k}{r^{\xi_k^-}}+\f{1}{\sqrt{\nu^2+4(k^2+i\mu k)}}\left(r^{\xi_k^+}\int_r^{\infty} s^{-\xi_k^++1} F_k(s)\mathrm{~d} s+r^{\xi_k^-}\int_1^r s^{-\xi_k^-+1}F_k(s)\mathrm{~d} s\right)\nonumber\\
    :=&{\bar{w}_k}{r^{\xi_k^-}}+\f{1}{\sqrt{\nu^2+4(k^2+i\mu k)}}h_{k,F}(r). \label{swk}
\end{align}

For $k\neq 0$, we intend to obtain  the constant $\bar{w}_k$ by applying the boundary condition of $\bl{\fv}$. After the $k$ mode of the vorticity mode is given, then from \eqref{streamvor}$_{1,3,4}$ we see that the $k$ mode of the stream function is given by
\be\label{Sphi}
\phi_k(r)={\bar{\phi}_k}{r^{-|k|}}+\f{r^{|k|}}{2|k|}\int_r^{\infty} s^{-|k|+1} w_k(s)\mathrm{d} s+\f{r^{-|k|}}{2|k|}\int_1^r s^{|k|+1} w_k(s)\mathrm{d} s\,,\q\text{for}\q k\in\mathbb{Z}-\{0\}\,.
\ee
Recall \eqref{stream1} and the boundary condition \eqref{PoNSvLF}$_4$, we derives
\be\label{Bf}
g_{r,k}=ik\phi_k(1),\q\text{and}\q g_{\th,k}=-\p_r\phi_k(1).
\ee
Substituting \eqref{Bf} in \eqref{Sphi}, we derive that
\[
\left\{
\begin{split}
\f{g_{r,k}}{ik}&=\bar{\phi}_k+\f{1}{2|k|}\int_1^{\infty}s^{-|k|+1} w_k(s)\mathrm{d} s\,,\\[1mm]
-g_{\th,k}&=-|k|\bar{\phi}_k+\f{1}{2}\int_1^\i s^{-|k|+1}w_k(s)\mathrm{d} s\,,
\end{split}
\right.
\]
and thus
\be\label{Sphik}
\left\{
\begin{split}
&\bar{\phi}_k=-\f{i}{2k}g_{r,k}+\f{1}{2|k|}g_{\th,k}\,,\\[1mm]
&\int_1^\i s^{-|k|+1}w_k(s)\md s=-ig_{r,k}\mathrm{sgn}\,k-g_{\th,k}\,,
\end{split}\quad \text { for }\q k \in \mathbb{Z}-\{0\}\,.
\right.
\ee
 Inserting \eqref{swk} into \eqref{Sphik}$_2$, one can deduce
\be\label{sbw}
\bar{w}_k=\left(g_{\th,k}+ig_{r,k}\mathrm{sgn}\,k+G_{k,F}\right)(2-|k|+\xi_k^{-})\,,
\ee
where
\be\label{GKF}
\begin{split}
G_{k,F}&:=\f{1}{\sqrt{\nu^2+4(k^2+i\mu k)}}\int_1^\i r^{-|k|+1}{h_{k,F}(r)}\mathrm{~d} r\,.\\
\end{split}
\ee
Therefore, the boundary condition for $\phi_k$ and $w_k$ is given by
\be\label{boundaryc}
\left\{
\begin{split}
&\bar{\phi}_k=\f{i}{2k}g_{r,k}+\f{1}{2|k|}g_{\th,k}\,,\\[1mm]
&\bar{w}_k=\left(g_{\th,k}-ig_{r,k}\mathrm{sgn}\,k+G_{k,F}\right)(2-|k|+\xi_k^{-})\,,
\end{split}\quad \text { for }\q k \in \mathbb{Z}-\{0\}\,.
\right.
\ee
Remembering \eqref{stream1}, for each $k\in\mathbb{Z}-\{0\}$, we can obtain $\fv_{r,k}(r)$ and $\fv_{\th,k}(r)$ by using \eqref{Sphi} and \eqref{boundaryc}$_1$ that:
\begin{small}
\be\label{EV}
\left\{\begin{split}
\fv_{r,k}(r)&=\f{1}{2}\left(g_{r,k}+i g_{\th,k}\mathrm{sgn\,}k\right)r^{-|k|-1}+\f{i\,\mathrm{sgn}k}{2}\left(r^{-|k|-1}\int_1^rs^{|k|+1}w_k(s)\mathrm{~d}s+r^{|k|-1}\int_r^\i s^{-|k|+1}w_k(s)\mathrm{~d}s\right)\,;\\[2mm]
\fv_{\th,k}(r)&=\f{1}{2}\left(g_{\th,k}-ig_{r,k}\mathrm{sgn\,}k\right)r^{-|k|-1}+\f{1}{2}\left(r^{-|k|-1}\int_1^rs^{|k|+1}w_k(s)\mathrm{~d}s-r^{|k|-1}\int_r^\i s^{-|k|+1}w_k(s)\mathrm{~d}s\right)\,.
\end{split}\right.
\ee
\end{small}


Now we estimate the $k-$mode solution in \eqref{EV} by $\bl{f}_k$. Set
\be\label{lamset}
\la=3^+<1-{\mathrm{Re}}\,\xi_1^{-}.
\ee
Due to the definition of $F_k$ given in \eqref{streamvor}$_2$, we first notice that
\be\label{GK1}
\begin{split}
\int_r^{\infty} s^{-\xi_k^++1} F_k(s)\mathrm{~d} s&=\int_r^{\infty} s^{-\xi_k^+} (sf_{\th,k}(s))'\mathrm{~d} s-ik\int_r^{\infty} s^{-\xi_k^+} f_{r,k}(s)\mathrm{~d} s\\
&=\xi_k^{+}\int_r^{\infty} s^{-\xi_k^+} f_{\th,k}(s)\mathrm{~d} s-r^{-\xi_k^++1}f_{\th,k}(r)-ik\int_r^{\infty} s^{-\xi_k^+} f_{r,k}(s)\mathrm{~d} s\,.
\end{split}
\ee
Meanwhile, it holds that
\be\label{GK2}
\begin{split}
\int_1^{r} s^{-\xi_k^-+1} F_k(s)\mathrm{~d} s&=\int_1^{r} s^{-\xi_k^-} (sf_{\th,k}(s))'\mathrm{~d} s-ik\int_1^{r} s^{-\xi_k^-} f_{r,k}(s)\mathrm{~d} s\\
&=\xi_k^{-}\int_1^{r} s^{-\xi_k^-} f_{\th,k}(s)\mathrm{~d} s+r^{-\xi_k^-+1}f_{\th,k}(r)-f_{\th,k}(1)-ik\int_1^{r} s^{-\xi_k^+} f_{r,k}(s)\mathrm{~d} s\,.
\end{split}
\ee

Combining \eqref{GK1} and \eqref{GK2}, one has
\be\label{HKF}
\begin{split}
h_{k,F}(r)=&\xi_k^{+}r^{\xi_k^+}\int_r^{\infty} s^{-\xi_k^+} f_{\th,k}(s)\mathrm{~d} s+\xi_k^{-}r^{\xi_k^-}\int_1^{r} s^{-\xi_k^-} f_{\th,k}(s)\mathrm{~d} s-f_{\th,k}(1)r^{\xi_k^-}\\
&-ik\left(r^{\xi_k^+}\int_r^{\infty} s^{-\xi_k^+} f_{r,k}(s)\mathrm{~d} s+r^{\xi_k^-}\int_1^{r} s^{-\xi_k^-} f_{r,k}(s)\mathrm{~d} s\right).
\end{split}
\ee
By denoting $
\|\bl{f}_k\|_{L^\i_{\la}}:=\|f_{\th,k}\|_{L^\i_{\la}}+\|f_{r,k}\|_{L^\i_{\la}}\,.
$ and using \eqref{lamset}, we can obtain that
\be\label{gkF}
\begin{split}
|h_{k,F}(r)|&\leq |f_{\th,k}(1)|\,r^{\mathrm{Re\,}\xi_k^-}+\|\bl{f}_k\|_{L^\i_{\la}}\left(\left(|\xi_k^+|+|k|\right)r^{\mathrm{Re\,}\xi_k^+}\int_r^\i s^{-\mathrm{Re\,}\xi_k^+-\la}\mathrm{~d} s+\left(|\xi_k^-|+|k|\right)r^{\mathrm{Re\,}\xi_k^-}\int_1^r s^{-\mathrm{Re\,}\xi_k^--\la}\mathrm{~d} s\right)\\
&= |f_{\th,k}(1)|\,r^{\mathrm{Re\,}\xi_k^-}+\|\bl{f}_k\|_{L^\i_{\la}}\left(\f{|\xi_k^+|+|k|}{\mathrm{Re\,}\xi_k^++\la-1}+\f{|\xi_k^-|+|k|}{-\mathrm{Re\,}\xi_k^--\la+1}\right)
r^{-\la+1}-\|\bl{f}_k\|_{L^\i_{\la}}\f{\left(|\xi_k^-|+|k|\right)}{-\mathrm{Re\,}\xi_k^--\la+1}r^{\mathrm{Re\,}\xi_k^-}\\
&\leq C\|\bl{f}_k\|_{L^\i_{\la}}\left(\f{|\xi_k^+|+|k|}{\mathrm{Re\,}\xi_k^++\la-1}+\f{|\xi_k^-|+|k|}{-\mathrm{Re\,}\xi_k^--\la+1}\right)r^{-\la+1}\,.
\end{split}
\ee
Substituting \eqref{gkF} in \eqref{GKF}, one deduces that
\[
\begin{split}
|G_{k,F}|\leq&{|\nu^2+4(k^2+i\mu k)|^{-1/2}}\|\bl{f}_k\|_{L^\i_{\la}}\left(\f{|\xi_k^+|+|k|}{\mathrm{Re\,}\xi_k^++\la-1}+\f{|\xi_k^-|+|k|}{-\mathrm{Re\,}\xi_k^--\la+1}\right)\int_1^\i r^{-|k|-\la+2}\mathrm{~d} r\\[1mm]
=&{|\nu^2+4(k^2+i\mu k)|^{-1/2}}(|k|+\la-1)^{-1}\left(\f{|\xi_k^+|+|k|}{\mathrm{Re\,}\xi_k^++\la-1}+\f{|\xi_k^-|+|k|}{-\mathrm{Re\,}\xi_k^--\la+1}\right)\|\bl{f}_k\|_{L^\i_{\la}}\,.
\end{split}
\]
For the convenience, we denote
\[
\left\{\begin{split}
&a_{k,\mu,\nu}:=\left|2-|k|+\xi_k^-\right|\,,\\[1mm]
&b_{k,\mu,\nu}:={|\nu^2+4(k^2+i\mu k)|^{-1/2}}(|k|+\la-1)^{-1}\left(\f{|\xi_k^+|+|k|}{\mathrm{Re\,}\xi_k^++\la-1}+\f{|\xi_k^-|+|k|}{-\mathrm{Re\,}\xi_k^--\la+1}\right)\,.
\end{split}
\right.
\]
From \eqref{xifirst} and \eqref{xisecond}, direct calculation implies that there exists a constant $C_{\nu,\mu}$ such that
\bes
|a_{k,\mu,\nu}|\leq C_{\mu,\nu}|k|,\q |b_{k,\mu,\nu}|\leq C_{\mu,\nu} |k|^{-2}.
\ees
Thus one finds that \eqref{sbw} infers
\be\label{wbk}
\begin{split}
|\bar{w}_k|\leq&\,|a_{k,\mu,\nu}|\left(|g_{\th,k}|+|g_{r,k}|+|G_{k,F}|\right)\leq C_{\mu,\nu}|k||\bl{g}_{k}|+C_{\mu,\nu}|k|^{-1}\|\bl{f}_k\|_{L^\i_{\la}}.
\end{split}
\ee
Now we are on the way to estimate the linear vorticity on $k-$mode. By \eqref{swk}, \eqref{gkF} and \eqref{wbk}, noticing that $\mathrm{Re\,}\xi_k^-+\la-1<0$ for any $k\neq0$, one derives
\be\label{evrt-1}
\begin{split}
|r^{\la-1}w_k(r)|\leq&|\bar{w}_k|r^{\mathrm{Re\,}\xi_k^-+\la-1}+{|\nu^2+4(k^2+i\mu k)|^{-1/2}}|h_{k,F}(r)|r^{\la-1}\\[2mm]
\leq&C_{\mu,\nu}|k||\bl{g}_{k}|+C_{\mu,\nu}|k|^{-1}\|\bl{f}_k\|_{L^\i_{\la}}.\\
\end{split}
\ee
By taking derivative on \eqref{swk}, we see that
\be\label{wpk}
w'_k(r)=\bar{w}_k \xi^{-}_k r^{\xi^-_k-1}+\f{1}{\sqrt{\nu^2+4(k^2+i\mu k)}}h'_{k,F}(r).
\ee
From \eqref{HKF}, it follows that
\be\label{hpfk}
\begin{split}
h'_{k,F}(r) =&\left(\xi_k^{+}\right)^2 r^{\xi_k^{+}-1} \int_r^{\infty} s^{-\xi_k^{+}} f_{\theta, k}(s) \md s+\left(\xi_k^{-}\right)^2 r^{\xi_k^{-}-1} \int_1^r s^{-\xi_k^{-}} f_{\theta, k}(s) \md s \\
& -\xi_k^{-} f_{\theta, k}(1) r^{\xi_k^{-}-1}-i k \xi_k^{+} r^{\xi_k^{+}-1} \int_{r}^{\infty} s^{-\xi_k^{+}} f_{r, k}(s) \md s-i k \xi_k^{-} r^{-\xi_k^--1}\int_1^r s^{-\xi_k^{-}} f_{r, k}(s) \md s\,.
\end{split}
\ee
Similarly as we achieve \eqref{gkF}, we derive from \eqref{hpfk} that
\[
|h_{k,F}(r)|\leq C_{\mu,\nu}|k| r^{-\la}\|\bl{f}_k\|_{L^\i_\la}\,.
\]
Thus we obtain from \eqref{wpk} that
\be\label{evrt0}
|r^{\la}w'_k(r)|\leq C_{\mu,\nu}k^2|\bl{g}_{k}|+C_{\mu,\nu}\|\bl{f}_k\|_{L^\i_{\la}}\,.
\ee
This gives the $L^\i_{\la-1}$ bound of $w_k$ and $L^\i_{\la}$ bound of $w'_k$. Recalling \eqref{EV}, one deduces
\be\label{evrt1}
\begin{split}
|r^{\la-2}\fv_{j,k}(r)|\leq&\f{1}{2}\|w_k\|_{L^\i_{\la-1}}\left(r^{-|k|-1}\int_1^rs^{|k|+2-\la}\mathrm{~d}s+r^{|k|-1}\int_r^\i s^{-|k|+2-\la}\mathrm{~d}s\right)+|g_{j,k}|r^{\la-|k|-3}\\[2mm]
\leq&\f{1}{2}\left(\f{1}{|k|+3-\la}+\f{1}{|k|-3+\la}\right)\|w_k\|_{L^\i_{\la-1}}+|g_{j,k}|r^{\la-|k|-3}\,,
\end{split}
\ee
for $j\in\{\th\,,\,r\}$. This indicates
\[
\begin{split}
\|\bl{\fv}_k\|_{L^\i_{\la-2}}\leq&\left(\frac{|k| a_{k, \mu, \nu}}{k^2-(\lambda-3)^2}+1\right)|\bl{g}_k|+\frac{|k| b_{k, \mu, \nu}}{k^2-(\lambda-3)^2}\left(a_{k, \mu, \nu}+|k|+\lambda-1\right)\|\bl{f}_k\|_{L_\lambda^{\infty}}\,\\[1mm]
                         \leq&C_{\mu,\nu}\left(|\bl{g}_k|+ |k|^{-2}\|\bl{f}_k\|_{L_\lambda^{\infty}}\right)\,.  \label{evrt2}
\end{split}
\]
Using the divergence-free condition and \eqref{stream1}, we see that
\[
\left\{
\begin{split}
&\fv'_{\theta, k}=-\frac{\fv_{\theta, k}}{r}+\frac{i k}{r} \fv_{r, k}+w_k\,; \\
&\fv'_{r, k}=-\frac{\fv_{r, k}}{r}-\frac{i k}{r} \fv_{\theta, k}\,,\end{split}\right.
\]
Also, by taking $r$ derivatives of the above equations,  and using \eqref{evrt-1}, \eqref{evrt0} and \eqref{evrt1}, we can conclude that
\bes
k^2\|\bl{\fv}_{k}\|_{L^\i_{\la-2}}+|k|\cd\|\bl{\fv}'_{k}\|_{L^\i_{\la-1}}+\|\bl{\fv}''_{k}\|_{L^\i_{\la}}\leq C_{\mu,\nu}\left(k^2|\bl{g}_k|+\|\bl{f}_k\|_{L_\lambda^{\infty}}\right)\,,\q\text{for all}\q k\in\mathbb{Z}-\{0\}\,.
\ees

\qed

\section{The nonlinear solution}\label{sec4}
We prove Theorem \ref{Main} in this section. First we consider the case that $\nu<-2$. Based on the derived linear estimates in the previous section, we solve the nonlinear problem in the Banach space $\mathcal{B}_\la$ given in \eqref{BLa} by applying the Banach fixed-point theorem, where $3<\la<\min\{\la_0,\la_1\}$. Given $\bar{\bl{v}}=\bar{v}_r\bl{e_r}+\bar{v}_\th\bl{e_\th}\in {\mathcal{B}}_\la$, consider the following linear system
\begin{small}
\be\label{PoNSvN}
\left\{\begin{split}
&-\left(\p_r^2+\f{1-\nu}{r}\p_r+\f{1}{r^2}\p_\th^2-\f{1}{r^2}\right)v_r+\f{\mu}{r^2}\p_\th v_r+\f{2}{r^2}\p_\th v_\th-\f{2\mu}{r^2}v_\th+\p_r\pi=\bar{f}_r\,,\\
& -\left(\p_r^2+\f{1-\nu}{r}\p_r+\f{1}{r^2}\p_\th^2-\f{1+\nu}{r^2}\right)v_\th+\f{\mu}{r^2}\p_\th v_\th-\f{2}{r^2}\p_\th v_r+\f{1}{r}\p_\th \pi=\bar{f}_\th\,,\\
&\p_\th v_\th+\p_r(rv_r)=0\,,\\
&v_\th\big|_{r=1}= g_\th,\q v_r\big|_{r=1}={g}_r, \q \bl{v}\big|_{r\rightarrow+\i}=0\,,
\end{split}\right.
\ee
\end{small}
where
\be\label{RHS}
\left\{
\begin{split}
&\bar{f}_r:=-\left(\bar{v}_r\p_r+\f{\bar{v}_\th}{r}\p_\th\right)\bar{v}_r+\f{\bar{v}_\th^2}{r}+f_r\,;\\[1mm]
&\bar{f}_\th:=-\left(\bar{v}_r\p_r+\f{\bar{v}_\th}{r}\p_\th\right)\bar{v}_\th-\f{\bar{v}_r\bar{v}_\th}{r}+f_\th\,.
\end{split}\right.
\ee
Since $\la>3$, the Young inequality for convolution:
\be\label{YC}
\sum_{n\in\mathbb{Z}} \big|\sum_{k\in\mathbb{Z}}a_kb_{n-k}\big|\leq \big(\sum_{k\in\mathbb{Z}}|a_k|\big)\big(\sum_{k\in\mathbb{Z}}|b_k|\big)
\ee
and \eqref{RHS}$_1$ indicate that
\be\label{nonlinear6}
\begin{split}
\|\bar{f}_{r}\|_{\mathcal{E}_\la}\leq & \left(\sum_{k \in \mathbb{Z}}\left\|\bar{v}_{r, k}\right\|_{L^{\infty}_{\lambda-2}}\right)\left(\sum_{k \in \mathbb{Z}}\left\|\bar{v}_{r, k}^{\prime}\right\|_{L^{\infty}_{\la-1}}\right)+\left(\sum_{k \in \mathbb{Z}}\left\|{\bar{v}}_{\theta, k}\right\|_{L_{\lambda-2}^{\infty}}\right) \left(\sum_{k \in \mathbb{Z}}|k| \cdot\left\|\bar{v}_{r, k}\right\|_{L^{\infty}_{\lambda-2}}\right)\\
&+\left(\sum_{k \in \mathbb{Z}}\left\|{\bar{v}}_{\theta, k}\right\|_{L_{\lambda-2}^{\infty}}\right)^2+\sum_{k \in \mathbb{Z}}\left\|f_{r, k}\right\|_{L^{\infty}_\lambda}\\
\leq&C\|\bar{\bl{v}}\|_{{B}_\la}^2+\|{f}_{r}\|_{\mathcal{E}_\la}\,.
\end{split}
\ee
Similarly one derives the same estimate for $\bar{f}_{\th}$ :
\be\label{nonlinear6.5}
\|\bar{f}_{\th}\|_{\mathcal{E}_\la}\leq C\|\bar{\bl{v}}\|_{{B}_\la}^2+\|{f}_{\th}\|_{\mathcal{E}_\la}\,.
\ee

On the other hand, for the case that $\nu\geq-2$, we solve the related nonlinear problem in the Banach space $\bar{\mathcal{B}}_\la$ given in \eqref{BLab}. For any
\[
\bar{\bl{v}}=\bar{v}_r\bl{e_r}+\bar{v}_\th\bl{e_\th}=\bar{v}_r\bl{e_r}+\left(\widetilde{\bar{v}}_\th+\f{\bar{\sigma}}{r}\right)\bl{e_\th}\in\bar{{\mathcal{B}}}_\la
\]
with $\bar{v}_r\bl{e_r}+\widetilde{\bar{v}}_\th\bl{e_\th}\in{\mathcal{B}}_\la$. Redefine
\be\label{RHS1}
\left\{
\begin{split}
&\bar{f}_r:=-\left(\bar{v}_r\p_r+\f{{\bar{v}}_\th}{r}\p_\th\right)\bar{v}_r+\f{{\bar{v}}_\th^2}{r}-\f{\bar{\sigma}^2}{r^3}+f_r\,;\\[1mm]
&\bar{f}_\th:=-\left(\bar{v}_r\p_r+\f{{\bar{v}}_\th}{r}\p_\th\right){\bar{v}}_\th-\f{\bar{v}_r{\bar{v}}_\th}{r}+f_\th\,.
\end{split}\right.
\ee
Here  the extra term $\f{\bar{\sigma}^2}{r^3}$ in \eqref{RHS1}$_1$ can be absorbed in the pressure $\pi$. Since $\la>3$, the Young inequality for convolution \eqref{YC} and \eqref{RHS}$_1$ indicate that
\be\label{nonlinear611}
\begin{split}
\|\bar{f}_{r}\|_{\mathcal{E}_\la}\leq & \left(\sum_{k \in \mathbb{Z}}\left\|\bar{v}_{r, k}\right\|_{L^{\infty}_{\lambda-2}}\right)\left(\sum_{k \in \mathbb{Z}}\left\|\bar{v}_{r, k}^{\prime}\right\|_{L^{\infty}_{\la-1}}\right)+\left(\sum_{k \in \mathbb{Z}}\left\|{\widetilde{\bar{v}}}_{\theta, k}\right\|_{L_{\lambda-2}^{\infty}}+\left|\bar{\sigma}\right|\right) \left(\sum_{k \in \mathbb{Z}}|k| \cdot\left\|\bar{v}_{r, k}\right\|_{L^{\infty}_{\lambda-2}}\right)\\
&+\sum_{k \in \mathbb{Z}}\left\|{\bar{v}}_{\theta, k}\right\|_{L_{\lambda-2}^{\infty}}\left(\sum_{k \in \mathbb{Z}}\left\|{\widetilde{\bar{v}}}_{\theta, k}\right\|_{L^{\infty}_{\la-2}}+2\left|\bar{\sigma}\right|\right)+\sum_{k \in \mathbb{Z}}\left\|f_{r, k}\right\|_{L^{\infty}_\lambda}\\
\leq&C\|\bar{\bl{v}}\|_{\bar{B}_\la}^2+\|{f}_{r}\|_{\mathcal{E}_\la}\,.
\end{split}
\ee
Similarly one derives the same estimate for $\bar{f}_{\th}$ :
\be\label{nonlinear6.511}
\|\bar{f}_{\th}\|_{\mathcal{E}_\la}\leq C\|\bar{\bl{v}}\|_{\bar{B}_\la}^2+\|{f}_{\th}\|_{\mathcal{E}_\la}\,.
\ee

In the following, we consider the both cases of $\nu<-2$ and $\nu\geq-2$ simultaneously by denoting
\[
\widetilde{\mathcal{B}}_\la=\left\{
\begin{split}
\mathcal{B}_\la,\q\text{for }\q\nu<-2\,;\\
\bar{\mathcal{B}}_\la,\q\text{for }\q\nu\geq-2\,.\\
\end{split}
\right.
\]
Combining \eqref{nonlinear6}, \eqref{nonlinear6.5}, \eqref{nonlinear611}, \eqref{nonlinear6.511},  and linear estimates in the previous section (Lemma \ref{lemzeromode} and Lemma \ref{lemnonzeromode}), we see there exists a unique $\bl{v}\in\widetilde{\mathcal{B}}_\la$ that solves \eqref{PoNSvN}, and it satisfies
\be\label{nonlinear9}
\|\bl{v}\|_{\widetilde{\mathcal{B}}_\la}\leq C\left(\|\bl{f}\|_{\mathcal{E}_\la}+\|\bl{g}\|_{\mathcal{V}}+\|\bar{\bl{v}}\|_{\widetilde{\mathcal{B}}_\la}^2\right)\,.
\ee
Therefore, given any $\bar{\bl{v}}\in\widetilde{\mathcal{B}_\la}$ with $\|\bar{\bl{v}}\|_{\widetilde{\mathcal{B}}_\la}<2C\e$, and $(\bl{f},\bl{g})\in\mathcal{E}_\la\times\mathcal{V}$ that $\|\bl{f}\|_{\mathcal{E}_\la}+\|\bl{g}\|_{\mathcal{V}}<\e$, where $\e<1/(4C^2)$, estimate \eqref{nonlinear9} indicates
\[
\|\bl{v}\|_{\widetilde{\mathcal{B}}_\la}\leq C\left(\e+4C^2\e^2\right)<2C\e\,.
\]
This indicates the solution map
\[
\mathcal{L}\q:\q\bar{\bl{v}}\q\to\q\bl{v}
\]
arise from \eqref{PoNSvN} maps the ball $\{\bl{u}\in\widetilde{\mathcal{B}}_\la\,:\,\|\bl{u}\|_{\widetilde{\mathcal{B}_\la}}<2C\e\}$ to itself. On the other hand, for any $\bar{\bl{v}}_1$, $\bar{\bl{v}}_2\in \{\bl{u}\in\widetilde{\mathcal{B}}_\la\,:\,\|\bl{u}\|_{\widetilde{\mathcal{B}}_\la}<2C\e\}$, we denote
\[
\begin{split}
\bl{h}:=&-\left\{\left(\bar{v}_{1,r}\p_r+\f{\bar{v}_{1,\th}}{r}\p_\th\right)\bar{v}_{1,r}-\f{\bar{v}_{1,\th}^2}{r}-\left(\bar{v}_{2,r}\p_r+\f{\bar{v}_{2,\th}}{r}\p_\th\right)\bar{v}_{2,r}+\f{\bar{v}_{2,\th}^2}{r}+\f{\left(\bar{\sigma}_1^2-\bar{\sigma}_2^2\right)}{r^3}\chi_{\nu\geq-2}\right\}\bl{e_r}\\
&-\left\{\left(\bar{v}_{1,r}\p_r+\f{\bar{v}_{1,\th}}{r}\p_\th\right)\bar{v}_{1,\th}+\f{\bar{v}_{1,r}\bar{v}_{1,\th}}{r}-\left(\bar{v}_{2,r}\p_r+\f{\bar{v}_{2,\th}}{r}\p_\th\right)\bar{v}_{2,\th}-\f{\bar{v}_{2,r}
\bar{v}_{2,\th}}{r}\right\}\bl{e_\th}\,.
\end{split}
\]
Clearly, the same estimate as obtain \eqref{nonlinear9} indicates that
\[
\begin{split}
\|\mathcal{L}\bar{\bl{v}}_1-\mathcal{L}\bar{\bl{v}}_2\|_{\widetilde{\mathcal{B}}_\la}\leq C\|\bl{h}\|_{\mathcal{E}_\la}\leq C\left(\|\bar{\bl{v}}_1\|_{\widetilde{\mathcal{B}}_\la}+\|\bar{\bl{v}}_2\|_{\widetilde{\mathcal{B}}_\la}\right)\|\bar{\bl{v}}_1-\bar{\bl{v}}_2\|_{\widetilde{\mathcal{B}}_\la}\leq 4C^2\e\|\bar{\bl{v}}_1-\bar{\bl{v}}_2\|_{\widetilde{\mathcal{B}}_\la}\,,
\end{split}
\]
which shows $\mathcal{L}$ is a contract mapping since $\e<(4C^2)^{-1}$. The Banach fixed point theorem shows there exists \emph{a unique solution} of the equation
\[
\mathcal{L}\bl{v}=\bl{v}
\]
in $\{\bl{u}\in\widetilde{\mathcal{B}_\la}\,:\,\|\bl{u}\|_{\widetilde{\mathcal{B}_\la}}<2C\e\}$. And it satisfies
\[
\|\bl{v}\|_{\widetilde{\mathcal{B}_\la}}\leq C\left(\|\bl{f}\|_{\mathcal{E}_\la}+\|\bl{g}\|_{\mathcal{V}}\right)\,.
\]
This completes the proof.

\section*{Acknowledgments}
\addcontentsline{toc}{section}{Acknowledgments}
\q The authors wish to thank Professors Zhengguang Guo and Matthieu Hillairet for helpful discussions on the non-uniqueness of the problem.

Z. Li is supported by the China Postdoctoral Science Foundation (No. 2024M763474) and the National Natural Science Foundation of China (No. 12001285). X. Pan is supported by the National Natural Science Foundation of China (No. 12031006, No. 12471222).

{\footnotesize

{\sc Z. Li: School of Mathematics and Statistics, Nanjing University of Information Science and Technology, Nanjing 210044, China, and Academy of Mathematics \& Systems Science, Chinese Academy of Sciences, Beijing 100190, China}

  {\it E-mail address:}  zijinli@nuist.edu.cn
\medskip

 {\sc X. Pan: School of Mathematics and Key Laboratory of MIIT, Nanjing University of Aeronautics and Astronautics, Nanjing 211106, China}

  {\it E-mail address:}  xinghong\_87@nuaa.edu.cn

}

\end{document}